\documentclass[12pt]{amsart}
\usepackage{latexsym,amssymb,amsmath,amscd,amsthm,amsxtra}
\usepackage{amsmath}
\usepackage{amsfonts}
\usepackage{amssymb}
\usepackage{graphicx}
\newtheorem{thm}{Theorem}[section]
\newtheorem{prop}[thm]{Proposition}
\newtheorem{lem}[thm]{Lemma}
\newtheorem{cor}[thm]{Corollary}
\newtheorem{dfn}[thm]{Definition}

\newtheorem{rmk}[thm]{Remark}

\newcommand{\complex}{{\mathbb C}}

\newcommand{\reals}{{\mathbb R}}

\newcommand{\calf}{{\mathcal F}}

\newcommand{\calp}{{\mathcal P}}

\newcommand{\frakg}{\mathfrak{g}}
\newcommand{\frakk}{\mathfrak{K}}

\newcommand{\bgl}{ {\mbox {\boldmath $gl$}}}

\newcommand{\C}{\mathbb C}

\newcommand{\sgn}{\text{sgn}}
\newcommand{\Sym}{\text{Sym}}
\newcommand{\gl}{\mathfrak{gl}}

\begin{document}
\title{Gelfand-Fuchs cohomology of invariant formal vector fields}
\author{Ilya Shapiro and Xiang Tang}
\date{}
\maketitle

\begin{abstract}
Let $\Gamma$ be a finite group acting linearly on a vector space
$V$. We compute the Lie algebra cohomology of the Lie algebra of
$\Gamma$-invariant formal vector fields on $V$. We use this
computation to define characteristic classes for foliations on
orbifolds.
\end{abstract}
\section{Introduction}
The space of smooth vector fields on a manifold $M$ is naturally
an infinite dimensional Lie algebra with the usual commutator
bracket. The continuous cohomology of this infinite dimensional
Lie algebra was studied by Gelfand-Fuchs (\cite{fuchs:book} and
references therein). It is now called the Gelfand-Fuchs cohomology
of $M$. The study of this cohomology leads to a better
understanding of the manifold $M$, in particular, it was proved by
Bott-Segal and Haefliger-Trauber that the Gelfand-Fuchs cohomology
of $M$ is a homotopy invariant of $M$ and is equal to the singular
cohomology of a space functorially constructed from $M$.

The Gelfand-Fuchs cohomology turned out to be a very useful tool
in the study of foliations.  Berstein-Rozenfeld and Gelfand-Fuchs
\cite{fuchs:book} applied it to the construction of the secondary
characteristic classes of foliations generalizing the
Godbillon-Vey class.  Connes and Connes-Moscovici,
\cite{connes:book} and \cite{connes-moscovici:index}, adopted the
Gelfand-Fuchs cohomology to the study of the transverse index
theory for foliations; a Hopf algebraic generalization of the
Gelfand-Fuchs cohomology was developed by them.

The connection between the Gelfand-Fuchs cohomology and the
algebraic index theory of formal deformation quantizations of
symplectic manifolds was developed by Nest-Tsygan
\cite{nest-tsygan:index} and Feigin-Felder-Shoikhet
\cite{ffs:index}. The present paper is motivated in part by the
second author's study of the algebraic index theory of orbifolds
\cite{ppt:index}. An equivariant version of the Gelfand-Fuchs
computation of the cohomology is a crucial step.

Let $V$ be a vector space over $\reals$ or $\complex$, and
$\Gamma$ a finite group acting linearly on $V$. Let
$W^\Gamma_\rho$ be the space of $\Gamma$-invariant formal vector
fields on $V$, where $\rho$ denotes the $\Gamma$ action on $V$. We
compute the Lie algebra cohomology of $W^\Gamma_\rho$ as well as
some of the relative cases.

In \cite{ppt:index} a very special case of the above question was
examined. Namely, some nonzero elements in
$H^\bullet(W^\Gamma_\rho, \gl(V)^\Gamma)$ are shown to exist when
$V$ is a complex vector space ($\gl(V)^\Gamma$ consists of
$\Gamma$-invariant linear transformations on $V$). In this paper,
we give a much more complete answer to the above question when $V$
is complex, or when $V$ is real and $\Gamma$ is cyclic (with a
technical assumption $\Gamma$ may be other than cyclic, but the
cyclic case is sufficient for an application to the characteristic
classes). The methods we use are similar to those of Gelfand-Fuchs
\cite{fuchs:book}, i.e. Hochschild-Serre spectral sequence and
invariant theory. However, we need to deal with a somewhat more
complicated algebras and $E_2$-terms. Our result compares the Lie
algebra cohomology to the cohomology of a certain truncated Weil
algebra as well as to the cohomology of a certain topological
space naturally associated to a classifying space. This
description of the result mirrors closely that of Gelfand-Fuchs.
We remark that our computation does depend on the field we work
with. The endomorphism algebra of an irreducible representation of
$\Gamma$ is $\complex$ when $V$ is complex, but may be $\reals$,
$\complex$, or $\mathbb{H}$ when $V$ is real. We assume $\Gamma$
to be cyclic to avoid the quaternionic case.

As an application, we follow Kontsevich's take on the
Gelfand-Fuchs' method \cite{fuchs:book, kontsevich} to define
characteristic classes for equivariant foliations. When a manifold
is equipped with an action of a finite group $\Gamma$, and a
foliation $\calf$ on $M$ is $\Gamma$-equivariant, $\calf$ descends
naturally to a foliation on the quotient space $M/\Gamma$. When
$M/\Gamma$ is not a manifold but an orbifold, we obtain a
foliation $\tilde{\calf}$ on the orbifold $M/\Gamma$. We prove
that a foliation on an orbifold $X=M/\Gamma$ actually induces a
foliation on the corresponding inertia orbifold
$\tilde{X}=\coprod_{\left<\gamma\right>\subset \Gamma}
M^\gamma/\Gamma^\gamma$, where $\left<\gamma\right>$ is the
conjugacy class\footnote{This notation is a possible source of
confusion as it is also used to denote the cyclic group generated
by $\gamma$, however its meaning is clear from the context.} of
$\gamma$ in $\Gamma$, $M^\gamma$ the fixed point submanifold of
$M$\footnote{Note that $M^\gamma$ may have quite different
connected components.}, and $\Gamma^\gamma$ consists of fixed
points of $\gamma$ in $\Gamma$ under the conjugation action.
Though we work always with $M/\Gamma$ our results hold generally
for arbitrary orbifolds $X$. Finally, we are able to define
characteristic classes for a foliation on an orbifold $X$ as
elements in the de Rham cohomology of its inertia orbifold
$\tilde{X}$.

The paper is arranged as follows. In Section 2, we compute the Lie
algebra cohomology of $W^\Gamma$. In Section 3, we prove that a
foliation on an orbifold defines a foliation on the corresponding
inertia orbifold; applying our computations of the Gelfand-Fuchs
cohomology, we define characteristic classes for foliations on
orbifolds.

\noindent\textbf{Acknowledgements.} Both authors are deeply
indebted to D.B. Fuchs for  encouragements and helpful advice, as
well as J. Kaminker, M. Pflaum, and H. Posthuma for useful
conversations. The research of the second author is partially
supported by NSF Grant 0604552.

\section{Cohomology of invariant subalgebras of formal vector fields}
Let $V$ be a vector space over $\reals$ or $\complex$, and $W_{V}$
be the Lie algebra of formal vector fields on $V$. Suppose we have
an action $\rho$ of a finite group $\Gamma$ on $V$, this induces
an action on $W_{V}$. We compute in this section the Lie algebra
cohomology of the Lie algebra $W_{\rho}^\Gamma$ of
$\Gamma$-invariant formal vector fields on $V$ when $V$ is a
complex vector space, or when $V$ is a real vector space and
$\Gamma$ is a cyclic group. We divide the computations into
several steps. In Section 2.1-2.3, we work out the case of a
complex vector space $V$, and in Section 2.4, we explain how to
extend our computation to the case of a real vector space $V$ and
a cyclic group $\Gamma$ action.

\subsection{Eigenvalues of the Euler field}
We point out that the $\Gamma$ action on $V$ can be made
unitary\footnote{Had $\Gamma$ not been a finite or more generally
compact group, the condition that the action be unitary would have
to be required.}. Accordingly, as unitary representations of
$\Gamma$ are completely reducible, $V$ can be split into
\[
V=V_0\oplus\bigoplus_{\alpha=1}^k m_\alpha W_{\alpha},
\]
where $V_0$ is the trivial $\Gamma$ representation, $W_\alpha$ is
irreducible, and $m_\alpha$ is the multiplicity of $W_\alpha$ in
$V$.

We introduce coordinates on $V$ as follows: $ x^i$ on  $V_0$,
$x^j_{\alpha, s}$ on $W_{\alpha, s}$ for all $1\leq \alpha\leq k,\
1\leq s\leq m_\alpha, 1\leq j\leq dim(W_\alpha)$. We consider the
following vector field
\[
\sum_{i, \alpha, s}x^i_{\alpha, s}\frac{\partial}{\partial
x^i_{\alpha, s}}
\]
that we denote by $X$. It is not difficult to see that $X$ is
$\Gamma$ invariant and therefore belongs to $W_\rho^\Gamma$.  We
need the following Lemma.

\begin{lem}
Let $\frakg$ be a Lie algebra and $X\in\frakg$ such that
$ad_\frakg X$ is diagonal with non-negative eigenvalues. Then the
inclusion of $\frakg^X$ into $\frakg$ induces and isomorphism on
Lie algebra cohomology.
\end{lem}
\begin{proof}
That $\frakg^X$ is a subalgebra follows from the Jacobi identity.
Since $X$ acts diagonally on $\frakg$ it does so also on
$\wedge^\bullet \frakg^*$ and each eigenspace is a subcomplex.
Since the action of $X$ on cohomology is trivial, only the
$0$-eigenspace subcomplex contributes non-trivially to the
cohomology.  By non-negativity of eigenvalues, the $0$-eigenspace
subcomplex is exactly $\wedge^\bullet (\frakg^X)^*$.
\end{proof}

\begin{cor}\label{cor:simplification}
Let $W_X=(W_{\rho}^\Gamma)^X$, then
$$H^\bullet(W_{\rho}^\Gamma)=H^\bullet(W_X)$$
\end{cor}

It is possible to give an explicit description of $W_X$ in terms
of the decomposition of $V$ into irreducibles. Namely, $W_X$ can
be identified with $W_{V_0}\ltimes \big(\bigoplus_{m_\alpha}
Poly(V_0)\otimes \frak{gl}_{m_\alpha}(\complex)\big)$, where
$W_{V_0}$ is the Lie algebra of formal vector fields on $V_0$, and
$Poly(V_0)$ is the algebra of polynomials on $V_0$. We remark that
$W_{V_0}$ acts on $Poly(V_0)$ naturally, and on $Poly(V_0)\otimes
\frak{gl}_{m_\alpha}(\complex)$ via the first factor. This action
defines the Lie bracket between $W_{V_0}$ and
$\bigoplus_{m_\alpha} Poly(V_0)\otimes
\frak{gl}_{m_\alpha}(\complex)$.  The bracket on $Poly(V_0)\otimes
\frak{gl}_{m_\alpha}(\complex)$ is extended $Poly(V_0)$-linearly
from that on $\frak{gl}_{m_\alpha}(\complex)$.

\subsection{Spectral sequence and invariant theory}
In this subsection, we will work with the case when there is only
one nonzero $m_\alpha$ in the representation $\rho$. Namely, we
compute the Lie algebra cohomology of $W_{V_0}\ltimes
Poly(V_0)\otimes \frak{gl}(W)$, where $V_0$, $W$ are some vector
spaces, and the Lie bracket between $W_{V_0}$ and
$Poly(V_0)\otimes \frak{gl}(W)$ is defined by the action of
$W_{V_0}$ on $Poly(V_0)$.

We consider the Lie subalgebra $\gl(V_0)\oplus\gl(W)$ in
$W_{V_0}\ltimes Poly(V_0)\otimes \frak{gl}(W)$, and use the
Hochschild-Serre spectral sequence to compute the Lie algebra
cohomology.

The $E_1$ term of this spectral sequence is as follows:
\begin{align*}
E^{p,q}_1&=H^q\left(\gl(V_0)\oplus\gl(W),\bigwedge^p\left(\dfrac{W_{V_0}\ltimes Poly(V_0)\otimes \frak{gl}(W)}{\gl(V_0)\oplus\gl(W)}\right)^*\right)\\
&=H^q\left(\gl(V_0)\oplus\gl(W),\bigwedge^p\left(\dfrac{\Sym
V_0^*\otimes V_0\oplus W^*\otimes W\otimes \Sym
V_0^*}{V_0^*\otimes V_0\oplus W^*\otimes W}\right)^*\right).
\end{align*}

Thus it is equal to
$H^\bullet(\gl(V_0)\oplus\gl(W),\C)\otimes\text{Inv}$, where Inv
stands for the following expression:
\[
\left(\bigwedge^\bullet V_0\otimes\bigwedge^\bullet (\Sym^{\geq
2}V_0^*\otimes V_0)\otimes\bigwedge^{\bullet}(W\otimes
W^*\otimes\Sym^{\geq 1}V_0^*)
\right)^{*\,\gl(V_0)\oplus\gl(W)}\mspace{-90mu}.
\]

Let us say a few words of justification for this manipulation.
Suppose that $V$, $V_i$ are some complex vector spaces. It is a
well known fact that $H^\bullet(\gl(V),M)=H^\bullet(\gl(V))\otimes
M^{\gl(V)}$ provided that $M$ is a tensor module, i.e. is a
submodule of $\bigotimes V\otimes\bigotimes V^*$.  The idea of the
proof is to introduce an inner product on $V$ which is used to
define the Casimir element $\Delta$ of the universal enveloping
algebra.  The Casimir element (belonging to the center and being
self adjoint) gives a decomposition of $\bigotimes
V\otimes\bigotimes V^*$ into
$\text{ker}\Delta\oplus\text{im}\Delta$ as a $\gl(V)$-module. It
is then apparent that the module of $\gl(V)$-invariants (which is
exactly $\text{ker}\Delta$) is isolated among these submodules as
the only submodule with the trivial infinitesimal character. Thus
it is the only one contributing non-trivially to the cohomology
(as is clear from the $Ext$ interpretation). These considerations
then easily apply also to $M$ itself.  Using the above ideas it is
immediate that
$H^\bullet(\bigoplus\gl(V_i),M)=H^\bullet(\bigoplus\gl(V_i))\otimes
M^{\bigoplus\gl(V_i)}$ provided that $M\subset
\bigotimes(\bigotimes V_i\otimes\bigotimes V_i^*)$, for $i\in I$.

At this point we need to describe the structure of the invariants
Inv as an algebra.

\begin{lem}
\begin{align*}
&\left(\bigwedge^\bullet V_0\otimes\bigwedge^\bullet (\Sym^{\geq
2}V_0^*\otimes V_0)\otimes\bigwedge^{\bullet}(W\otimes
W^*\otimes\Sym^{\geq 1}V_0^*) \right)^{*\,\gl(V_0)\oplus\gl(W)}\\
&\mspace{60mu}=\sum_{r,s}\left(\bigwedge^{r+s}V_0\otimes\bigwedge^r(\Sym^2
V_0^*\otimes V_0)\otimes\bigwedge^s (W\otimes W^*\otimes
V_0^*)\right)^{*\,\gl(V_0)\oplus\gl(W)}\mspace{-90mu}.
\end{align*}

\end{lem}
\begin{proof}
This is proved in Proposition 5.2 \cite{ppt:index}. Here we
mention the ingredients. First of all by $\gl(V_0)$ invariance,
elements in $V_0$ have to be paired with elements in $V_0^*$. That
is any homogeneous component of an element will have equal numbers
of $V_0$ and $V_0^*$. Next, we observe that one cannot pair
$\wedge^k V_0$ with $\Sym^k V_0^*$ in a $\gl(V_0)$ invariant way
for $k>1$.  This is sufficient to prove the Lemma.

\end{proof}

Using the $\gl(V_0)\oplus\gl(W)$-invariance (c. f.
\cite{fuchs:book}), we observe that
$$\left(\bigwedge^{r+s}V_0\otimes\bigwedge^r(\Sym^2 V_0^*\otimes
V_0)\otimes\bigwedge^s (W\otimes W^*\otimes
V_0^*)\right)^{*\,\gl(V_0)\oplus\gl(W)}$$ is spanned by
$\psi_{\sigma,\gamma}$ with $\sigma \in \Sigma_r,\ \gamma\in
\Sigma_s$ and
\begin{align*}
\psi_{\sigma,\gamma}&(v_1,...,v_{r+s},\varphi_1,v_1^{'},...,\varphi_r,v_r^{'},...,
w_1,w_1^*,v_1^*,...,w_s,w_s^*,v_s^*)\\
&=\mspace{-36.0mu}\displaystyle\sum_{\alpha,\beta,\delta\in
\Sigma_{r+s}\times\Sigma_r\times\Sigma_s}\mspace{-36.0mu}\sgn(\alpha)\sgn(\beta)\sgn(\delta)\prod_{i=1}^r
\varphi_{\beta(i)}(v^{'}_{\beta\sigma(i)},v_{\alpha(i)})
\prod_{j=1}^s
w^*_{\delta(j)}(w_{\delta\gamma(j)})v^*_{\delta(j)}(v_{\alpha(r+j)})
\end{align*}
where $v_i, v_i'\in V_0$, $\varphi_i\in \Sym^2 V_0^*$, $w_j\in W$,
$w_j^*\in W^*$, and $v_k^*\in V_0^*$.

We will use the notation
$\psi_{\sigma,\gamma}((v_i)_1^{r+s},(\varphi_i,v_i^{'})_1^r,(
w_i,w_i^*,v_i^*)_1^s)$ to stand for the above pairing.  The way to
generate such a formula is to write down an element without any
symmetry conditions first, i.e.
$$\prod_{i=1}^r \varphi_i(v^{'}_{\sigma(i)},v_{\alpha(i)})
\prod_{j=1}^s w^*_j(w_{\gamma(j)})v^*_j(v_{\alpha(r+j)})$$ and
then anti-symmetrize it.  Re-indexing the above summations, we get
another form of the same element $\psi_{\sigma, \gamma}$ that is
evidently constant on the conjugacy classes of $\sigma$ and
$\gamma$:
$$\displaystyle\sum_{\alpha,\beta,\delta\in
\Sigma_{r+s}\times\Sigma_r\times\Sigma_s}\mspace{-36.0mu}\sgn(\alpha)\prod_{i=1}^r
\varphi_i(v^{'}_{\beta\sigma\beta^{-1}(i)},v_{\alpha(i)})
\prod_{j=1}^s
w^*_j(w_{\delta\gamma\delta^{-1}(j)})v^*_j(v_{\alpha(r+j)})$$

In order to better understand $\psi_{\sigma,\gamma}$, we define
$\Phi_\sigma\in \left(\bigwedge^{r}V_0\otimes\bigwedge^r(\Sym^2
V_0^*\otimes V_0)\right)^{*\,\gl(V_0)}$ and
$\Psi_\gamma\in\left(\bigwedge^{s}V_0\otimes\bigwedge^s (W\otimes
W^*\otimes V_0^*)\right)^{*\,\gl(V_0)\oplus\gl(W)}$ by
$$\Phi_\sigma((v_i)^r_1,(\varphi_i,v_i^{'})_1^r)=\displaystyle
\sum_{\nu,\beta\in\Sigma_r}\sgn(\nu)\prod_{i=1}^r\varphi_i
(v^{'}_{\beta\sigma\beta^{-1}(i)},v_{\nu(i)})$$
$$\Psi_\gamma((v_i)_1^s,(w_i,w_i^*,v_i^*)_1^s)=\displaystyle
\sum_{\eta,\omega\in\Sigma_s}\sgn(\eta)\prod_{j=1}^s
w^*_j(w_{\omega\gamma\omega^{-1}(j)})v^*_j(v_{\eta(j)})$$

We then have the following:
\begin{lem}
$$\Phi_\sigma\cdot\Psi_\gamma=\psi_{\sigma,\gamma}$$
\end{lem}
\begin{proof}
\begin{align*}
&\Phi_\sigma\cdot\Psi_\gamma((v_i)^{r+s}_1,(\varphi_i,v_i^{'})^r_1,(w_i,w_i^*,v_i^*)_1^s)\\
&=\frac{1}{r!s!}\displaystyle\sum_{\alpha\in\Sigma_{r+s}}\sgn(\alpha)\Phi_\sigma((v_{\alpha(i)})_1^r,(\varphi_i,v_i^{'})_1^r)\Psi_\gamma((v_{\alpha(r+i)})_1^s,(w_i,w_i^*,v_i^*)_1^s)\\
&=\frac{1}{r!s!}\displaystyle\mspace{-10.0mu}\sum_{\alpha\in\Sigma_{r+s}}\mspace{-10.0mu}\sgn(\alpha)\mspace{-10.0mu}\displaystyle\sum_{\nu,\beta\in\Sigma_r}\mspace{-10.0mu}\sgn(\nu)\prod_{i=1}^r\varphi_i(v^{'}_{\beta\sigma\beta^{-1}(i)},v_{\alpha\nu(i)})\mspace{-10.0mu}\displaystyle\sum_{\eta,\omega\in\Sigma_s}\mspace{-10.0mu}\sgn(\eta)\prod_{j=1}^sw^*_j(w_{\omega\gamma\omega^{-1}(j)})v^*_j(v_{\alpha(r+\eta(j))})\\
&=\frac{1}{r!s!}\displaystyle\sum_{\nu,\eta}\sgn(\nu\eta)\displaystyle\sum_{\alpha,\beta,\omega}\sgn(\alpha)\prod_{i=1}^r\varphi_i(v^{'}_{\beta\sigma\beta^{-1}(i)},v_{\alpha\nu(i)})\prod_{j=1}^sw^*_j(w_{\omega\gamma\omega^{-1}(j)})v^*_j(v_{\alpha(r+\eta(j))})\\
\intertext{after re-indexing}
&=\frac{1}{r!s!}\displaystyle\sum_{\nu,\eta} \displaystyle
\sum_{\alpha,\beta,\omega}\sgn(\alpha) \prod_{i=1}^r\varphi_i
(v^{'}_{\beta\sigma\beta^{-1}(i)},v_{\alpha(i)})\prod_{j=1}^s
w^*_j(w_{\omega\gamma\omega^{-1}(j)})v^*_j(v_{\alpha(r+j)})\\
&=\psi_{\sigma,\gamma}((v_i)^{r+s}_1,(\varphi_i,v_i^{'})^r_1,(w_i,w_i^*,v_i^*)_1^s).
\end{align*}
\end{proof}

The element $\Phi_\sigma$ obviously depends only on the conjugacy
class of $\sigma$.  Furthermore, if $\sigma\in\Sigma_r$ can be
decomposed as $\sigma_1\coprod\sigma_2$ with
$\sigma_i\in\Sigma_{r_i}$ and $r_1+r_2=r$, then
$\Phi_\sigma=\Phi_{\sigma_1}\cdot\Phi_{\sigma_2}$. This was proved
in Theorem 2.1.4 \cite{fuchs:book}. Thus
$\left(\bigwedge^{\bullet}V_0\otimes\bigwedge^{\bullet}(\Sym^2
V^*_0\otimes V_0)\right)^{*\,\gl(V_0)}$ is generated by
$\Phi_i\in\left(\bigwedge^{i}V_0\otimes\bigwedge^{i}(\Sym^2
V^*_0\otimes V_0)\right)^{*\,\gl(V_0)}$ for
$i=1,...,\text{dim}V_0$, given by the elementary cycles of length
$i$.

To deal with $\Psi_\gamma$, we consider the map:
$$\widetilde{}\quad:\left(\Sym^s(W\otimes W^*)\right)^{*\,\gl(W)}\rightarrow \left(\bigwedge^s V_0\otimes \bigwedge^s (W\otimes W^*\otimes V_0^*)\right)^{*\,\gl(V_0)\oplus\gl(W)}$$
defined by
$$\widetilde{\Psi}((v_i)_1^s,(w_i,w_i^*,v_i^*)_1^s)=\Omega_s((v_i)_1^s,(v_i^*)_1^s)\Psi((w_i,w_i^*)_1^s)$$
where
$$\Omega_s((v_i)_1^s,(v_i^*)_1^s)=\displaystyle\sum_{\alpha\in\Sigma_s}\sgn(\alpha)\prod_{i=1}^s
v_i^*(v_{\alpha(i)})$$ is the canonical element of $(\bigwedge^s
V_0\otimes \bigwedge^s V_0^*)^{*\,\gl(V_0)}$.

\begin{lem}
The map\quad$\widetilde{}$\quad is a homomorphism.
\end{lem}
\begin{proof}
Let $\Psi^s\in(\Sym^s(W\otimes W^*))^{*\,\gl(W)}$ and
$\Psi^{s'}\in(\Sym^{s'}(W\otimes W^*))^{*\,\gl(W)}$. Then
\begin{align*}
\widetilde{\Psi^s}&\cdot\widetilde{\Psi^{s'}}((v_i)_1^{s+s'},(w_i,w_i^*,v_i^*)_1^{s+s'})\\
&=\frac{1}{(s!)^2(s'!)^2}\sum_{\alpha,\beta\in\Sigma_{s+s'}}\sgn(\alpha\beta)\widetilde{\Psi^s}((v_{\alpha(i)})_1^s,(w_{\beta(i)},w^*_{\beta(i)},v^*_{\beta(i)})_1^s)\\
&\qquad\cdot\widetilde{\Psi^{s'}}((v_{\alpha(s+i)})_1^{s'},(w_{\beta(s+i)},w^*_{\beta(s+i)},v^*_{\beta(s+i)})_1^{s'})\\
&=\frac{1}{s!s'!}\sum_{\beta}\sgn(\beta)\Psi^s((w_{\beta(i)},w^*_{\beta(i)})_1^s)\Psi^{s'}((w_{\beta(s+i)},w^*_{\beta(s+i)})_1^{s'})\\
&\qquad\cdot\frac{1}{s!s'!}\sum_{\alpha}\sgn(\alpha)\Omega_s((v_{\alpha(i)})_1^s,(v^*_{\beta(i)})_1^s)\Omega_{s'}((v_{\alpha(s+i)})_1^{s'},(v^*_{\beta(s+i)})_1^{s'})\\
&=\Omega_{s+s'}((v_i)_1^{s+s'},(v^*_i)_1^{s+s'})\frac{1}{s!s'!}\sum_{\beta}\Psi^s((w_{\beta(i)},w^*_{\beta(i)})_1^s)\Psi^{s'}((w_{\beta(s+i)},w^*_{\beta(s+i)})_1^{s'})\\
&=\widetilde{\Psi^s\cdot\Psi^{s'}}((v_i)_1^{s+s'},(w_i,w_i^*,v_i^*)_1^{s+s'})
\end{align*}
\end{proof}

The homomorphism\quad$\widetilde{}$\quad obviously factors through
the quotient of the symmetric algebra by the ideal of functions of
degree greater than dim$(V_0)$.  The induced map from the quotient
is than an isomorphism as it is easy to construct an inverse.
Namely, consider the evaluation map on any partial basis of $V_0$
(this map is choice independent):

$$\text{ev}:\left(\bigwedge^s
V_0\otimes \bigwedge^s (W\otimes W^*\otimes
V^*_0)\right)^{*\,\gl(V_0)\oplus\gl(W)}\longrightarrow
\left(\Sym^s(W\otimes W^*)\right)^{*\,\gl(W)}$$
$$\Psi\mapsto\Psi((e_i)_1^s,(-,-,e_i^*)_1^s)$$ Then
clearly $\text{ev}\circ\widetilde{}=\text{Id}$, and to check that
$\,\,\widetilde{}\circ\text{ev}=\text{Id}$ it is sufficient to
consider $\Psi_\gamma$.  Recall that
$\Psi_\gamma((v_i)_1^s,(w_i,w_i^*,v_i^*)_1^s)=\mspace{-10.0mu}\displaystyle
\sum_{\eta,\omega\in\Sigma_s}\!\!\!\sgn(\eta)\prod_{j=1}^s
w^*_j(w_{\omega\gamma\omega^{-1}(j)})v^*_j(v_{\eta(j)})$ and
consider the element $\Psi\in\left(\Sym^s(W\otimes
W^*)\right)^{*\,\gl(W)}$ defined by
$\Psi((w_i,w_i^*)_1^s)=\displaystyle
\sum_{\omega\in\Sigma_s}\prod_{j=1}^s
w^*_j(w_{\omega\gamma\omega^{-1}(j)})$.  Then
$\Psi_\gamma=\widetilde{\Psi}$ and this is enough.

Since the structure of $\left(\Sym^\bullet(W\otimes
W^*)\right)^{*\,\gl(W)}$ is well known\footnote{The generators are
$\Psi_i\in\left(\Sym^i(W\otimes W^*)\right)^{*\,\gl(W)}$ for
$i=1,...,\text{dim}(W)$ corresponding to the elementary cycles of
length $i$, and there are no relations.}, we obtain the structure
of $\left(\bigwedge^\bullet V_0\otimes \bigwedge^\bullet (W\otimes
W^*\otimes V_0^*)\right)^{*\,\gl(V_0)\oplus\gl(W)}$ as an algebra.
Namely it has generators $\widetilde{\Psi}_i$ and relations
$\prod\widetilde{\Psi}_{i_k}=0$ if $\sum i_k>\text{dim}(V_0)$.
Note that might mean that some of the generators are themselves
$0$. We are now ready for:

\begin{lem}
\label{lem:Inv} Inv is generated by $\Phi_i,\ i=1,\cdots,dim(V_0)$
of degree $2i$ and $\widetilde{\Psi}_j,\ j=1,\cdots,dim(W)$ of
degree $2j$ subject to the relation that the total degree can not
exceed $2dim(V_0)$.
\end{lem}
\begin{proof}
The only thing that remains is to show that there are no other
relations.  Suppose that
$\displaystyle\sum_{I,J}a_{I,J}\Phi_I\widetilde{\Psi}_J=0$, where
$I$ and $J$ are multi-indices and
$|I|+|J|\leq\text{dim}(V)$.\footnote{By $|I|$ we mean $\sum i_k$
where $I=\{i_1,...,i_n\}$.} We may assume that for all $I$ and
$J$, $|I|=r$ and $|J|=s$, thus we may relabel the sum to
$\displaystyle\sum_{[\sigma],J}a_{[\sigma],J}\Phi_{\sigma}\widetilde{\Psi}_J=0$
with $[\sigma]$ running over the conjugacy classes of $\Sigma_r$.
Choosing a partial basis $(e_i)_1^{r+s}$ of $V$, and a
$\tau\in\Sigma_r$, set $\varphi_{e_i}=(e^*_i)^2$ and consider the
following calculation:
\begin{align*}
&\Phi_{\sigma}\widetilde{\Psi}_J((e_i)_1^{r+s},(\varphi_{e_i},e_{\tau^{-1}(i)})_1^r,(w_i,w_i^*,e^*_{r+i})_1^s)\\
&=\frac{1}{r!s!}\sum_{\alpha\in\Sigma_{r+s}}\sgn(\alpha)\Phi_\sigma((e_{\alpha(i)})_1^r,(\varphi_{e_i},e_{\tau^{-1}(i)})_1^r)\Omega_s((e_{\alpha(r+i)})_1^s,(e^*_{r+i})_1^s)\Psi_J((w_i,w_i^*)_1^s).\\
\intertext{Note that
$\Omega_s((e_{\alpha(r+i)})_1^s,(e^*_{r+i})_1^s)=0$ unless
$\alpha=\alpha_1\coprod\alpha_2$ with $\alpha_1\in\Sigma_r$ and
$\alpha_2\in\Sigma_s$, in which case it equals $\sgn(\alpha_2)$,
thus}
&=\Psi_J((w_i,w_i^*)_1^s)\frac{1}{r!}\sum_{\beta\in\Sigma_r}\sgn(\beta)\Phi_\sigma((e_{\beta(i)})_1^r,(\varphi_{e_i},e_{\tau^{-1}(i)})_1^r)\\
&=\Psi_J((w_i,w_i^*)_1^s)\Phi_\sigma((e_i)_1^r,(\varphi_{e_i},e_{\tau^{-1}(i)})_1^r)
\end{align*}
Recalling the definition we see that
$$\Phi_\sigma((e_i)_1^r,(\varphi_{e_i},e_{\tau^{-1}(i)})_1^r)=\displaystyle
\sum_{\nu,\beta\in\Sigma_r}\sgn(\nu)\prod_{i=1}^r\varphi_{e_i}
(e_{\tau^{-1}\beta\sigma\beta^{-1}(i)},e_{\nu(i)})=
\begin{cases}
|\text{Stab}(\sigma)| & \text{if $\sigma\sim\tau$}\\
0 & \text{else.}
\end{cases}$$  Therefore evaluating
$\displaystyle\sum_{[\sigma],J}a_{[\sigma],J}\Phi_{\sigma}\widetilde{\Psi}_J$
at
$(e_i)_1^{r+s},(\varphi_{e_i},e_{\tau^{-1}(i)})_1^r,(-,-,e^*_{r+i})_1^s$
we obtain that $\displaystyle\sum_J a_{[\tau],J}
|\text{Stab}(\tau)| \Psi_J=0$ in $(\Sym^\cdot (W\otimes
W^*))^{*\,\gl(W)}$, thus $a_{[\tau],J}=0$ for all $\tau$ and $J$.
\end{proof}
\begin{rmk}\label{rmk:relgl}
The above Lemma is actually a description of
$H^\bullet(W^\Gamma_n,GL^\Gamma_n)$ as an algebra.  See also Prop.
\ref{prop:rel-complex}
\end{rmk}

\subsection{Weil algebras and classifying spaces}
In the next step, we relate the cohomology of the Lie algebra
$W_{\rho}^\Gamma$ to the cohomology of a truncated Weil algebra as
well as to the cohomology of a certain topological space.

We recall the general definition of a Weil algebra. Let $\frakg$
be a Lie algebra. Define the Weil algebra $W(\frakg)$ to be
$\wedge^\bullet \frakg^*\otimes S^\bullet \frakg^*$. Introduce a
grading on $W(\frakg)$ by assigning elements from $\wedge^i
\frakg^*\otimes S^j \frakg^*$ degree $i+2j$. The algebra
$W(\frakg)$ is filtered by a decreasing filtration
\[
W(\frakg)=F^0W(\frakg)\supset F^1W(\frakg)=F^2W(\frakg)\supset
F^3W(\frakg)=F^4W(\frakg)\supset\cdots,
\]
with $F^pW(\frakg)=\bigoplus\limits_{2j\geq p} \wedge^\bullet
\frakg^*\otimes S^j \frakg^*$.

A differential $d_W$ on $W(\frakg)$ can be introduced as follows:
\[
\begin{split}
&d_W\psi(g_1\wedge \cdots\wedge g_i\otimes h_1\otimes \cdots h_j)=\sum_{t=1}^j\psi(g_1\wedge \cdots\wedge g_i\wedge h_t\otimes h_1\otimes \cdots \hat{h}_t\otimes \cdots \otimes h_j)\\
&+\sum_{s=1}^{i}\sum_{t=1}^j(-1)^{s-1}\psi([g_s,h_t]\wedge g_1\wedge \cdots\hat{g}_s\cdots \wedge g_i\otimes h_1\otimes \cdots\hat{h}_t\cdots \otimes h_j)\\
&+\sum_{1\leq s_1<s_2\leq
i}(-1)^{s_1+s_2-1}\psi([g_{s_1},g_{s_2}]\wedge g_1\wedge
\cdots\hat{g}_{s_1} \cdots \hat{g}_{s_2} \cdots\wedge g_i\otimes
h_1\otimes \cdots\otimes h_j),
\end{split}
\]
and $(W(\frakg),d_W, F^\bullet)$ forms a filtered differential
graded algebra. One computes the $E_2$ term of the spectral
sequence associated to the filtration:
\[
E_2^{p,q}=\left\{\begin{array}{ll}0,&p\ \text{is odd},\\
H^q(\frakg; S^{p/2}\frakg^*),&p\ \text{is even}.\end{array}\right.
\]
We note that $(W(\frakg), d)$ is acyclic.

The Weil algebra $W(\frak{g})$ has a lifting property which we
will use below. Suppose that $(C^\bullet, \delta)$ is a
differential graded algebra. Then any linear map $f$ from
$\frakg^*$ to $C^1$ of $(C^\bullet, \delta)$ can be lifted to a
homogeneous multiplicative homomorphism $F$ from $(W(\frakg),
d_W)$ to $(C^\bullet,\delta)$, which agrees with $f$ on $\frakg
^*$. To define such a map $F$, it is sufficient to specify the
image of the generators:
\[
\begin{split}
F(g^*\otimes 1)=f(g^*),\ \ \ \ \ \ & F(1\otimes
g^*)=-\delta(f(g^*)).
\end{split}
\]

For the purpose of our computation, we need to consider the
truncated Weil algebra
$W(\frakg)_{2n}=W(\frakg)/F_{2n+1}W(\frakg)$. The filtration
$F^\bullet$ descends to a filtration on $W(\frakg)_{2n}$, and the
$E_2$-term of the spectral sequence associated to the inherited
filtration is
\[
E_2^{p,q}=\left\{\begin{array}{ll}H^q(\frakg; S^{p/2}\frakg^*),&p\
\text{is even and }p\leq 2n,\\
0,&\text{otherwise}.\end{array}\right.
\]

We consider the special case of $\frakg=\frak{gl}(V_0)\oplus
\frak{gl}(W)$. To compute the cohomology of the truncated Weil
algebra of $\gl(V_0)\oplus \gl(W)$, we need to understand the Lie
algebra cohomology of $\frak{gl}(V_0)\oplus \frak{gl}(W)$ with
coefficient in $S^\bullet(\frak{gl}(V_0)\oplus \frak{gl}(W))^*$.
The cohomology $H^\bullet(\frak{gl}(V_0)\oplus \frak{gl}(W)$; $
S^{\bullet}(\frak{gl}(V_0)\oplus\frak{gl}(W))^*)$ is isomorphic to
\[
H^\bullet(\frak{gl}(V_0);S^\bullet \frak{gl}(V_0)^*)\otimes
H^\bullet(\frak{gl}(W); S^\bullet \frak{gl}(W)^*).
\]
By the Casimir element argument, as is explained in Section 2.2,
we see that
\[
\begin{split}
H^\bullet(\frak{gl}(V_0);S^\bullet \frak{gl}(V_0)^*)&=H^\bullet(\frak{gl}(V_0))\otimes (S^\bullet\frak{gl}(V_0)^*)^{\frak{gl}(V_0)},\\
H^\bullet(\frak{gl}(W); S^\bullet
\frak{gl}(W)^*)&=H^\bullet(\frak{gl}(W))\otimes
(S^\bullet\frak{gl}(W)^*)^{\frak{gl}(W)}.
\end{split}
\]
Furthermore by \cite{fuchs:book}[Theorem 2.1.5],
$(S^\bullet\frak{gl}(V_0)^*)^{\frak{gl}(V_0)}$ is generated by
$\xi_i$, $i=1,\cdots,$ $dim(V_0)$, with $deg(\xi_i)=2i$ and
\[
\xi_i(g_1, \cdots, g_i)=\sum_{\sigma\in
\Sigma_i}Tr(g_{\sigma(1)}\cdots g_{\sigma(i)}),
\]
for $g_s\in \frak{gl}(V_0)$. Here we view an element $g$ in
$\frak{gl}(V_0)$ as a linear endomorphism on $V_0$, and $Tr$ is
the trace functional on linear endomorphisms. The same holds for
$(S^\bullet\frak{gl}(W)^*)^{\frak{gl}(W)}$. Summarizing, we
observe that
\[
\begin{split}
&H^\bullet(\frak{gl}(V_0)\oplus \frak{gl}(W);
S^\bullet(\frak{gl}(V_0)\oplus \frak{gl}(W))^*)\\
=&H^\bullet(\gl(V_0))\otimes H^\bullet(\gl(W))\otimes
\complex[\xi_1, \cdots, \xi_{dim(V_0)}]\otimes \complex[\eta_1,
\cdots, \eta_{dim(W)}].
\end{split}
\]

Let us compute the Lie algebra cohomology of $W_{V_0}\ltimes
Poly(V_0)\otimes \gl(W)$ using the Weil algebra. Consider the
``natural" projection $pr$ from $W_{V_0}\ltimes
Poly(V_0)\otimes\frak{gl}(W)$ to $\frak{gl}(V_0)\oplus
\frak{gl}(W)$. The dual of this map defines a linear map from
$(\frak{gl}(V_0)\oplus \frak{gl}(W))^*$ to $(W_{V_0}\ltimes
Poly(V_0)\otimes\frak{gl}(W))^*$. By the universal lifting
property of the Weil algebra we obtain a map
\[
\chi:W(\frak{gl}(V_0)\oplus \frak{gl}(W))\rightarrow
\wedge^\bullet(W_{V_0}\ltimes Poly(V_0)\otimes \frak{gl}(W))^*.
\]
We write out $\chi$ explicitly on the generators of
$W(\gl(V_0)\oplus\gl(W))$,
\begin{equation}\label{eq:chi}
\begin{split}
\chi(g^*\otimes 1)=g^*,\ \ \ \ \ \ &
\chi(h^*\otimes 1)=h^*,\\
\chi(1\otimes g^*)=-dg^*,\ \ \ \ \ & \chi(1\otimes h^*)=-dh^*,
\end{split}
\end{equation}
for $g\in\frak{gl}(V_0)$, $h\in \frak{gl}(W)$.

We notice that for $g\in \frak{gl}(V_0)$, $h\in \frak{gl}(W)$,
$dg^*\in V_0\wedge \Sym^2V_0^*\otimes V_0$ and $dh^*\in V_0\wedge
V^*_0\otimes W^*\otimes W$. It is easy to see that the $\chi$
image of $F^{2n+1}W(\frak{gl}(V_0)\oplus \frak{gl}(W))$ will be
contained in $\wedge^{\geq n+1}V_0\otimes \wedge^\bullet\cdots$.
When $n+1$ is greater than the dimension of $V_0$, then the above
is automatically zero. Therefore, the map $\chi$ factors through
the quotient $F^{2\dim(V_0)+1}W(\frak{gl}(V_0)\oplus
\frak{gl}(W))$. Furthermore, observing that both $\chi(1\otimes
g^*)=-dg^*$ and $\chi(1\otimes h^*)=-dh^*$ are contained in
\[
\wedge^2\big(\frac{\Sym V_0^*\otimes V_0\oplus W^*\otimes W\otimes
\Sym V_0^*}{V_0^*\otimes V_0\oplus W^*\oplus W}\big)^*,
\]
and both $\chi(g^*\otimes 1)=g^*$ and $\chi(h^*\otimes 1)=h^*$ are
contained in $\frak{gl}(V_0)\oplus \frak{gl}(W)$, we conclude that
the map $\chi$ is compatible with the filtrations on
$W(\frak{gl}(V_0)\oplus \frak{gl}(W))$ and
$\wedge^\bullet(W_{V_0}\ltimes Poly(V_0)\otimes \frak{gl}(W))^*$.
We remind the reader that the filtration on
$\wedge^\bullet(W_{V_0}\ltimes Poly(V_0)\otimes \frak{gl}(W))^*$
(its associated spectral sequence is the Hochschild-Serre spectral
sequence for the Lie subalgebra $\frak{gl}(V_0)\oplus
\frak{gl}(W)$) is defined by the powers of the wedge
\[
\wedge^\bullet\big(\frac{\Sym V_0^*\otimes V_0\oplus W^*\otimes
W\otimes \Sym V_0^*}{V_0^*\otimes V_0\oplus W^*\otimes W}\big)^*.
\]

In conclusion, we have constructed a filtration compatible
differential graded algebra morphism
\[
\chi:W(\frak{gl}(V_0)\oplus \frak{gl}(W))_{2\dim(V_0)}\rightarrow
\wedge^\bullet(W_{V_0}\ltimes Poly(V_0)\otimes \frak{gl}(W))^*,
\]
where $W(\frak{gl}(V_0)\oplus \frak{gl}(W))_{2\dim(V_0)}$ is the
truncated Weil algebra.

\begin{lem}
\label{lem:weil} The map $\chi$ defined above is a
quasi-isomorphism.
\end{lem}
\begin{proof}
Since $\chi$ is compatible with the filtrations, it is sufficient
to prove that $\chi$ is an isomorphism on the $E_2$ terms of the
spectral sequences associated to the filtrations. As we know that
$\chi$ is multiplicative, it is sufficient to show that $\chi$
defines isomorphisms on $E^{p,0}_2$ and $E^{0,q}_2$.

For $E^{0,q}_2$, we easily see from Equation (\ref{eq:chi}) that
$\chi$ is an identity map on $\wedge^\bullet \gl(V_0)\otimes
\wedge^\bullet\gl(W)$.

For $E^{p,0}_2$, we observe that $E^{\bullet,0}_2$ forms an
algebra. And as $\chi$ is multiplicative, this reduces to checking
the statement on the generators. Observe that $E_2^{\bullet,0}$ of
$W(\frak{gl}(V_0)\oplus \frak{gl}(W))$ is generated by $\xi_i$ and
$\eta_j$, $i=1,\cdots, \dim(V_0)$ and $j=1, \cdots \dim(W)$. We
compute $\chi(\xi_i)$ and $\chi(\eta_j)$.

\begin{enumerate}
\item The image $\chi(\xi_i)$ was already considered in the proof
of Theorem 2.2.4', \cite{fuchs:book}. Hence, we skip the
particulars of this part. The conclusion is that up to a non-$0$
constant the image of $\chi(\xi_i)$ is $\Phi_i$.

\item For $\chi(\eta_j)$, looking at Equation (\ref{eq:chi}), we
observe that $\chi(\eta_j)$ is in
\[
\wedge^jV_0\otimes \wedge^j (W^*\otimes W\otimes V^*),
\]
as $\chi(1\otimes h)$ is in $V_0\otimes W^*\otimes W\otimes
V_0^*$.

Recall that for an element $h\in \gl(W)^*$, $\chi(h)(v, w\otimes
w^*\otimes v')=h([v,w\otimes w^*\otimes v' ])=v'(v)h(w\otimes
w^*)$, for $v\in V_0$, $v'\in V^*$, $w\in W$, and $w^*\in W^*$. As
$\chi$ is an algebra homomorphism, we generalize the above
evaluation to $\chi(\eta_j)$, i.e.
\[
\begin{split}
&\chi(\eta_j)(v_1,\cdots, v_j, w_1^*, w_1', v_1^*,\cdots, w_j^*, w_j', v_j^* )\\
=&\prod_{\sigma\in \Sigma_j}\sgn
(\sigma)v_1^*(v_{\sigma(1)})\cdots
v_j^*(v_{\sigma(j)})\eta_j(w_1^*\otimes w_1',\cdots, w_j^*\otimes
w_j').
\end{split}
\]
We recall that $\eta_j$ is the element in $(\Sym^j
\gl(W)^*)^{\gl(W)}$ corresponding to the $j$-cycle in the
permutation group $\Sigma_j$. The second line of the above
equation agrees with the $\widetilde{\Psi}_j$'s evaluation on
$v_1,\cdots, v_j, w_1^*, w_1, v_1^*,\cdots, w_j^*, w_j, v_j^*$.
Therefore, we have that $\chi(\eta_j)=\widetilde{\Psi}_j$.
\end{enumerate}

With Lemma \ref{lem:Inv}, we conclude that $\chi$ is an
isomorphism on the $E_2$ terms of the spectral sequences and
therefore induces an isomorphism on the cohomologies.
\end{proof}

Having dealt with the special case of a single non-$0$ $m_\alpha$,
we can state the result in full generality.  The proof is
identical though more notationally intensive and so we omit it.

\begin{thm}
\label{thm:main}Let $\Gamma$ be a finite group acting linearly on
a complex vector space $V$. Suppose that $V$ decomposes as
$V=V_0\oplus\bigoplus_{\alpha=1}^km_\alpha W_\alpha$ into
irreducible representations of $\Gamma$, with $V_0=V^\Gamma$. Then
\[
H^\bullet(W_\rho^\Gamma)=H^\bullet(W(\gl(V_0)\oplus
\bigoplus_{\alpha=1}^k\gl_{m_\alpha}(\complex))_{2\dim(V_0)}).
\]
\end{thm}

We now wish to give a topological description of
$H^\bullet(W_\rho^\Gamma)$. Consider the principal bundle
$$\pi: E(GL(V_0)\times \prod_{\alpha=1}^k GL_{m_\alpha}(\complex))\to B(GL(V_0)\times \prod_{\alpha=1}^k GL_{m_\alpha}(\complex)).
$$
Denote by $B_\rho^\Gamma$ the $2\dim(V_0)$-th skeleton  of
$B(GL(V_0)\times \prod_{\alpha=1}^k GL_{m_\alpha}(\complex))$, and
by $X^\Gamma_{\rho}$ its $\pi$ preimage inside $E(GL(V_0)\times
\prod_{\alpha=1}^k GL_{m_\alpha}(\complex))$.  The following
Theorem is then almost an immediate corollary of Thm.
\ref{thm:main}

\begin{thm}\label{thm:main-top}
$H^\bullet(W_{\rho}^\Gamma)\cong H^\bullet(X^\Gamma_{\rho})$.
\end{thm}

\begin{proof}
The proof is identical to that of Theorem 2.2.4',
\cite{fuchs:book}. Here we outline the main ingredients. We
compare the spectral sequence associated to the filtration on
$W(\gl(V_0)\oplus\bigoplus_{\alpha=1}^k
\gl_{m_\alpha}(\complex))_{2dim(V_0)}$ with the spectral sequence
of the bundle $X_\rho^\Gamma\rightarrow B_\rho^\Gamma$. From the
previous computation, we have seen that $E_2$ terms of the two
spectral sequences agree as algebras. Furthermore, by the
acyclicity of the untruncated Weil algebra
$W(\gl(V_0)\oplus\bigoplus_\alpha\gl_{m_\alpha}(\complex))$ and
the cohomology of $E(GL(V_0)\times \prod_{\alpha=1}^k
GL_{m_\alpha}(\complex))$, we see that in the two spectral
sequences the exterior generators of the algebra $$
E^{0,\bullet}_2=H^\bullet(\gl(V_0)\oplus
\bigoplus_{\alpha=1}^k\gl_{m_\alpha}(\complex))=H^\bullet(GL(V_0)
\times\prod_{\alpha=1}^k
GL_{m_\alpha}(\complex))
$$
are transgressive and mapped to the multiplicative generators of
the algebra
$$
E^{\bullet,0}_2=(S^\bullet(\gl(V_0)\oplus\bigoplus_{\alpha=1}^k
\gl_{m_\alpha}(\complex))^*)^{\gl(V_0)\oplus\bigoplus_\alpha\gl_{m_\alpha}
(\complex)}_{2dim(V_0)}=H^\bullet(B_\rho^\Gamma)
$$
by the transgression. Therefor, we conclude that the limits of the
two spectral sequences agree.

\end{proof}

\begin{rmk}
When both $V$ and $W$ are 1 dimensional, the Lie algebra
cohomology of $W_V\ltimes Poly(V)\otimes \gl(W)$ can be computed
directly without using spectral sequence. Let $x$ be the variable
on $V$ and $y$ be the variable on $W$. Consider the vector field
$x\frac{\partial}{\partial x}\in \gl(V)$. We find that the adjoint
action of $x\frac{\partial}{\partial x}$ on $W_V\ltimes
Poly(V)\otimes \gl(W)$ and therefore on the cochain complex is
diagonal. Accordingly, the cohomology of $W_V\ltimes
Poly(V)\otimes \gl(W)$ is computed by the 0 eigenvectors in
\[
\bigwedge^\bullet \big(W_V\ltimes Poly(V)\otimes \gl(W)\big)^*.
\]
We can find these 0 eigenvectors explicitly and compute the Lie
algebra cohomology easily without using spectral sequence.
\end{rmk}

We now have the language to state the following Proposition:

\begin{prop}
\label{prop:rel-complex} $H^\bullet(W_{\rho}^\Gamma,
\frak{gl}(V_0)\oplus\bigoplus_{\alpha=1}^k\frak{gl}_{m_\alpha}(\complex))
=H^\bullet(B_\rho^\Gamma)$.
\end{prop}

\subsection{The real case}\label{sec:realcase}
In this subsection, we consider the real Lie algebra
$W_{\rho}^\Gamma$ of $\Gamma$ invariant formal vector fields on a
real vector space $V$, where $\rho$ denotes the $\Gamma$
representation.  Let us assume that $\Gamma$ is a finite cyclic
group as that is the case we need for defining characteristic
classes of foliation on an orbifold.

As $\Gamma$ is finite, the $\Gamma$ action on $V$ is completely
reducible. Suppose that $U$ is an irreducible component. There are
two possibilities: $U\otimes\C$ is irreducible, thus
$End_\Gamma(U\otimes\C)=\C$ and so $End_\Gamma(V)=\reals$ or
$U\otimes\C=W\oplus\overline{W}$ with $W$ irreducible in which
case $End_\Gamma(U\otimes\C)=\C^2$ since
$W\neq\overline{W}$\footnote{This is where one needs the
assumption that $\Gamma$ is cyclic.  Alternatively, we may simply
assume that $W\neq\overline{W}$.  This rules out the quaternions.}
(they correspond to different eigenvalues of the generator of
$\Gamma$), and so $End_\Gamma(V)=\C$. Therefore, we can write,
similar to the complex case:
\[
V=V_0\oplus m_{-1}W_{-1}\oplus\bigoplus_{\alpha=1}^k m_\alpha
W_\alpha,
\]
where $\Gamma$ acts on $V_0$ trivially, $W_{-1}$ is the one
dimensional representation with character -1, $W_\alpha$ is an
irreducible representation of $\Gamma$ with
$dim_\reals(W_\alpha)=2$, and $m_\alpha$ is the multiplicity of
the representation $W_\alpha$ in $V$.

Similar to the complex case in Section 2.1, there is a
$\Gamma$-invariant Euler vector field $X$ on $\reals^n$, which
acts on $W_V$ diagonally with nonnegative eigenvalues. As with
Corollary \ref{cor:simplification}, by looking at the eigenvalues
of $X$ we reduce the Lie algebra $W_{\rho}^\Gamma$ to $W_X$, where
$W_X$ consists of eigenvectors of $X$ in $W_\rho^\Gamma$ with zero
eigenvalue. We see that $W_X$ is somewhat different from the
complex case:
\[
W_X=W_{V_0}\ltimes Poly(V_0)\otimes \big(\gl_{m_{-1}}(\reals)
\oplus\bigoplus_{\alpha=1}^k \bgl_{m_\alpha}\big),
\]
where $\bgl_{m_\alpha}$ is the complex general Lie algebra
$\frak{gl}_{m_\alpha}(\complex)$ viewed as a real Lie algebra.

Extending results \ref{cor:simplification}-\ref{lem:weil}, we have
\begin{prop}
\label{prop:real-coh}There is a natural quasi-isomorphism $\chi$
from the truncated Weil algebra $$W\big(\frak{gl}(V_0)\oplus
\gl_{m_{-1}}(\reals)\oplus\bigoplus_{\alpha=1}^k \bgl _{m_\alpha}
\big)_{2dim(V_0)}$$ to the cochain complex of the real Lie algebra
$W_{\rho}^\Gamma$.
\end{prop}
\begin{proof}
The proof is a copy of that of Lemma \ref{lem:weil}. The map
$\chi$ is defined by the projection from $W_X$ to $\gl(V_0)\oplus
\gl_{m_{-1}}(\reals)\oplus \bigoplus_{\alpha=1}^k
\bgl_{m_{\alpha}}$. To prove that $\chi$ is a quasi-isomorphism,
we compare the $E_2$ terms of spectral sequences associated to the
filtration on the Weil algebra $W(\gl(V_0)\oplus
\gl_{m_{-1}}(\reals)\oplus\bigoplus_{\alpha=1}^k
\bgl_{m_{\alpha}})_{2dim(V_0)}$ and to the Lie subalgebra
$\gl(V_0)\oplus \gl_{m_{-1}}(\reals)\oplus \bigoplus_{\alpha=1}^k
\bgl_{m_{\alpha}}$ of $W_ X$. One change we need to make is to
replace $\gl_{m_\alpha}(\complex)$ by $\bgl_{m_\alpha}$, while the
algebra of invariant polynomials on $\bgl_{m_\alpha}$ is computed
in Proposition \ref{prop:bgl-symm}.
\end{proof}

Furthermore, we may choose a $\Gamma$-invariant metric on $V$ and
consider the Lie subalgebra
$\frak{o}(V_0)\oplus\frak{o}_{m_{-1}}(\reals)\oplus\bigoplus_{\alpha=1}^k
\frak{u}_{m_\alpha}$ inside $W_{\rho}^\Gamma$, where $\frak{o}$ is
the Lie algebra of orthogonal matrices, and $\frak{u}$ is the Lie
algebra of skew-hermitian matrices.

We now use the relative version of the Hochschild-Serre spectral
sequence to compute the relative cohomology
$H^\bullet(W_{\rho}^\Gamma,
\frak{o}(V_0)\oplus\frak{o}_{m_{-1}}(\reals)\oplus\bigoplus_{\alpha=1}^k\frak{u}_{m_\alpha})$.
We need the following Lemma:
\begin{lem}\label{lem:invariant}
Let $V_i$ be real vector spaces, $W_i$ complex vector spaces,
$\frak{o}_i$ and $\frak{u}_i$ reductive subalgebras of $\gl(V_i)$
and $\bgl(W_i)$ respectively, then
\begin{align*}
H^\bullet(\bigoplus\gl(V_i)&\oplus\bigoplus\bgl(W_i),
\bigoplus\frak{o}_i\oplus\bigoplus\frak{u}_i;M)\\
&=H^\bullet(\bigoplus\gl(V_i)\oplus\bigoplus\bgl(W_i),
\bigoplus\frak{o}_i\oplus\bigoplus\frak{u}_i)\otimes
M^{\bigoplus\gl(V_i)\oplus\bigoplus\bgl(W_i)}
\end{align*}
provided $M$ is a ``tensor" module.
\end{lem}

\begin{proof}
Observe that
$H^\bullet(\frak{g},\frak{k},M)=H^\bullet(\frak{g},\frak{k})\otimes
M^{\frak{g}}$ if its complexified version holds.  Note also that
$\bgl(W)\otimes\C=\gl_\C(W)\oplus\gl_\C(W)$; this decomposition is
given by the eigenspaces of the ``forgotten" multiplication by $i$
acting on $\bgl(W)$.  The meaning of a ``tensor" module is now
clear. We may use the $Ext$ interpretation of the relative
cohomology (valid for reductive subalgebras, see
\cite{borel-wallach:book} for example) to conclude that only a
module with a trivial infinitesimal character contributes and thus
the earlier discussion of cohomology of tensor modules over
$\bigoplus\gl(V_i)$ applies.
\end{proof}

\begin{prop}
\label{prop:relative-coh}Let
\[
W\big(\gl(V_0)\oplus
\gl_{m_{-1}}(\reals)\oplus\bigoplus_{\alpha=1}^k\bgl_{m_\alpha},
\frak{o}(V_0)\oplus\frak{o}_{m_{-1}}(\reals)\oplus
\bigoplus_{\alpha=1}^k\frak{u}_{m_\alpha}\big)_{2dim(V_0)}
\]
be the truncated (at degree $>2dim(V_0)$) relative Weil algebra.
The map $\chi$ descends to a quasi-isomorphism to the relative
cochain complex $C^\bullet(W_\rho^\Gamma, \frak{o}(V_0)\oplus
\frak{o}_{m_{-1}}(\reals)\oplus\bigoplus_{\alpha=1}^n\frak{u}_{m_\alpha})$.
\end{prop}
\begin{proof}
The proof consists of checking that the $\chi$ is a
quasi-isomorphism on the $E_2$ terms associated to the relative
spectral sequences, which is analogous to the proof of Lemma
\ref{lem:weil}.
\end{proof}

To obtain a topological description of the cohomology  we compute
the invariants of $\bgl_m$ as follows.
\begin{prop}
\label{prop:bgl-symm} There is an isomorphism of algebras
$$(S^\bullet\bgl_m)^{*\bgl_m}=\reals[x_s,y_s]_{s=1}^m$$ with
degree of $x_s$ and $y_s$ equal to $s$.

\end{prop}
\begin{proof}
As is explained in the proof of Lemma \ref{lem:invariant}, the
complexification of $\bgl_m$ is isomorphic to
$\gl_m{\complex}\oplus\gl_m(\complex)$. While
$(S^\bullet\gl_m(\complex))^{*\gl_m(\complex)}=\C[z_s]_{s=1}^m$
with $z_s$ of degree $s$ given by the symmetric trace, and so
$(S^\bullet \bgl_m\otimes \complex)^{*\bgl_m\otimes \complex}$ is
a complex polynomial algebra with $2m$ generators $z_s$ and $t_s$.
It is straight forward to trace through the isomorphisms to check
that $(S^\bullet \bgl_m)^{*\bgl_m}$ (as the subalgebra on which
conjugation acts trivially) is isomorphic to
$\reals[\frac{z_s+t_s}{2},\frac{z_s-t_s}{2i}]$.
\end{proof}

\begin{prop}\label{prop:bgl-u} There exist isomorphisms:
\begin{enumerate}
\item $H^\bullet(\bgl_m)=H^\bullet(U_m\times U_m)$, \item
$H^\bullet(\bgl_m, \frak{u}_m)=H^\bullet(U_m)$.
\end{enumerate}
\end{prop}
\begin{proof}
For statement (1), we notice that $\bgl_m\otimes
\complex=\gl_m(\complex)\oplus\gl_m(\complex)=(\frak{u}_m\oplus\frak{u}_m)\otimes
\complex$. Therefore,
$H^\bullet(\bgl_m)=H^\bullet(\frak{u}_m\oplus
\frak{u}_m)=H^\bullet(U_m\times U_m)$ because $U_m\times U_m$ is a
compact Lie group.

For statement (2), we observe that under the isomorphism
$\bgl_m\otimes \complex=\gl_m(\complex)\oplus\gl_m(\complex)$, the
complexification of $\frak{u}_m$ is identified with
$\gl_m(\complex)$ embedded diagonally into $\bgl_m\otimes
\complex$. Furthermore, the pair $(\gl_m(\complex)\oplus\gl_m(\complex),
\gl_m(\complex))$ is the complexification of
$(\frak{u}_m\oplus\frak{u}_m, \frak{u}_m)$ with $\frak{u}_m$
 embedded into $\frak{u}_m\oplus\frak{u}_m$ diagonally.
Therefore, $ H^\bullet(\bgl_m,
\frak{u}_m)=H^\bullet(\frak{u}_m\oplus\frak{u}_m,
\frak{u}_m)=H^\bullet(U_m\times U_m/U_m)=H^\bullet(U_m).$
\end{proof}

By Proposition \ref{prop:bgl-symm}, $(S^\bullet\bgl_m)^{*\bgl_m}$
is isomorphic to the cohomology ring of the classifying space
$B(U_m\times U_m)$, a polynomial ring with $2m$ generators. This
brings to mind the fibration $E(U_m\times U_m)\rightarrow
B(U_m\times U_m)$ with the fiber $U_m\times U_m$. Thus we obtain
the following topological description of the cohomology computed
in Proposition \ref{prop:real-coh} and \ref{prop:relative-coh}.
Consider the fibration
\begin{align*}
\pi:E(U_{dim(V_0)}\times U_{m_{-1}}\times&\prod_{\alpha=1}^k
(U_{m_\alpha}\times U_{m_\alpha}))\\
&\longrightarrow B(U_{dim(V_0)}\times
U_{m_{-1}}\times\prod_{\alpha=1}^k(U_{m_\alpha}\times
U_{m_\alpha})).
\end{align*}
Let $X_\rho^\Gamma$ be the $\pi$ preimage in $E(U_{dim(V_0)}\times
U_{m_{-1}}\times\prod_{\alpha=1}^k (U_{m_\alpha}\times
U_{m_\alpha}))$ of the $2dim(V_0)$-th skeleton $B_\rho^\Gamma$ of
$B(U_{dim(V_0)}\times
U_{m_{-1}}\times\prod_{\alpha=1}^k(U_{m_\alpha}\times
U_{m_\alpha}))$.  Note that $SO(V_0)\times
SO_{m_{-1}}(\reals)\times \prod_{\alpha=1}^k U_{m_\alpha}$ and
$O(V_0)\times O_{m_{-1}}(\reals)\times \prod_{\alpha=1}^k
U_{m_\alpha}$ act on the fibers of $X_\rho^\Gamma\rightarrow
B_\rho^\Gamma$ ($U_{m_{\alpha}}$'s are embedded diagonally in
$U_{m_\alpha}\times U_{m_\alpha}$).  We then have:

\begin{thm}
\label{thm:top-real}
$$H^\bullet(W_{\rho}^\Gamma)=H^\bullet(X_\rho^\Gamma);$$
$$H^\bullet(W_{\rho}^\Gamma,\gl(V_0)\oplus
\gl_{m_{-1}}(\reals)\oplus\bigoplus_{\alpha=1}^k\bgl_{m_\alpha})
=H^\bullet(B_\rho^\Gamma);$$
$$H^\bullet(W_{\rho}^\Gamma,\frak{o}(V_0)\oplus
\frak{o}_{m_{-1}}(\reals)\oplus
\bigoplus_{\alpha=1}^k\frak{u}_{m_\alpha})=H^\bullet(X_\rho^\Gamma/SO(V_0)\times
SO_{m_{-1}}(\reals)\times \prod_{\alpha=1}^k U_{m_\alpha});$$
$$H^\bullet(W_\rho^\Gamma, O(V_0)\times O_{m_{-1}}(\reals)\times
\prod_{\alpha=1}^k
U_{m_\alpha})=H^\bullet(X^\Gamma_\rho/O(V_0)\times
O_{m_{-1}}(\reals)\times \prod_{\alpha=1}^k U_{m_\alpha}).$$

\end{thm}
\begin{proof}
The proof of the first statement is similar to that of Theorem
\ref{thm:main-top}, ones uses Proposition \ref{prop:real-coh}. The
second statement was essentially demonstrated in the process of
proving \ref{thm:main-top}.  With Proposition
\ref{prop:relative-coh}, the last two statements follow by
comparing the cohomologies of the Weil algebras and the
topological spaces. Compare to \cite{fuchs:book} Thm. 2.2.6.

\end{proof}

\section{Characteristic Classes for foliations on orbifolds}
In this section, we apply the computations of the Lie algebra
cohomology above to define some characteristic classes for
foliations on orbifolds.
\subsection{Foliations on orbifolds}
Here, we introduce a notion of a foliation on an orbifold. Because
all of the following constructions and computations are local,
instead of general orbifolds, we will work with the global
quotient, i.e. $X=M/\Gamma$, where $\Gamma$ is a finite group
acting on a smooth manifold $M$. We consider a foliation $\calf$
on $M$, which is invariant under the $\Gamma$ action, and call it
a $\Gamma$ equivariant foliation on $M$.  We denote by $\calf$
both the foliation and the distribution defining the foliation.

\begin{thm}\label{thm:foliation}Let $\calf$ be a $\Gamma$ equivariant
foliation on $M$. Then for each $\gamma\in \Gamma$, $\calf$
restricts to a foliation $\calf^\gamma$ on $M^\gamma$, where
$M^\gamma$ is the $\gamma$ fixed point manifold.
\end{thm}
\begin{proof}
Note that $\calf^\gamma$ on $M^\gamma$ is given by
$\calf^\gamma_x=T_x M^\gamma\cap \calf_x$. We will prove that
$\calf^\gamma$ is indeed a foliation on $M^\gamma$, i.e.
\begin{enumerate}
\item $\calf^\gamma$ is of constant rank on each connected
component of $M^\gamma$; \item $\calf^\gamma$ is integrable.
\end{enumerate}
We check (1) and (2) separately.

For (1), we notice that for $x\in M^\gamma$, $\gamma$ acts on $T_x
M$ with $\calf_x$ a submodule, and $\calf^\gamma_x=(\calf_x)^{G}$.
Where $G=\left<\gamma\right>$, a finite cyclic group generated by
$\gamma$. Since in the representation space of $G$, the trivial
representation is an isolated point\footnote{This phrasing of the
proof is motivated by the remark that follows it.  To prove only
the statement of the theorem one can observe that the eigenvalues
of $\gamma$ are $n$-th roots of unity and so cannot change
continuously.  This in effect proves the discreteness of the
representation space.}, the dimension of $\calf_x^\gamma$ is
locally constant on $M^\gamma$. Therefore, $\calf^\gamma$ is of
constant rank on each component of $M^\gamma$.

To check (2), we observe that a pair of sections $\eta$, $\nu$ of
$\calf^\gamma$ can locally be extended to a pair $\tilde{\eta}$,
$\tilde{\nu}$ of $G$-invariant sections of $\calf$. Since $\calf$
is closed under the commutator bracket, and the $G$-action
commutes with the bracket, we see that
$[\tilde{\eta},\tilde{\nu}]$ is a $G$-invariant section of
$\calf$.  Thus $[\eta,\nu]=[\tilde{\eta},\tilde{\nu}]|_{M^\gamma}$
is a section of $\calf^\gamma$.
\end{proof}

\begin{rmk}
As is clear from the proof, we can extend Theorem
\ref{thm:foliation} to the following setting. Let $G$ act on a
manifold $M$ such that $M$ has a $G$-invariant hermitian
structure. Suppose that $\calf$ is a $G$-invariant foliation on
$M$. Then $\calf$ restricts to a foliation on the fixed point
manifold $M^G$ when $G$ satisfies property $T$.
\end{rmk}
\subsection{Characteristic classes}
In this subsection, we want to define some characteristic classes
for a foliation $\calf$ on an orbifold $X$. As foreshadowed by
Theorem \ref{thm:foliation}, our characteristic classes map takes
values in the cohomology group of the inertia orbifold
$H^\bullet(\tilde{X})$.

We review the theory of characteristic classes for foliations on a
manifold. Possible references for the standard approach to this
are \cite{fuchs:book}, \cite{kt:foliation}. Our approach differs
and owes much to the point of view of \cite{kontsevich}.  In
addition we develop the equivariant version of these methods. Let
$G$ be a Lie group with its Lie algebra denoted by $\frakg$.  Our
discussion takes place over $\reals$, i.e. all manifolds and Lie
algebras are real.  It readily modifies, where appropriate, to
$\C$ by considering almost complex manifolds and complex Lie
algebras.

Recall that a $\frakg$-structure on a smooth manifold $M$ is a
smooth 1-form $\omega$ on $M$ with values in $\frakg$, satisfying
the Maurer-Cartan equation:
\[
d\omega+\frac{1}{2}[\omega, \omega]=0.
\] The $\frakg$-structure on $M$ is equivalent to the structure of a
trivialized flat principal $G$-bundle $\calp=M\times G$ over $M$.
The flat connection on $\calp$ defines a $\frakg$-structure on
$M$.

Given a $\frakg$-structure $\omega$ on a manifold $M$, we can
define a characteristic classes map $\chi_\omega$ by
\[
\begin{array}{cccc}
\chi_\omega:& \wedge^\bullet \frakg^* &\rightarrow &\Omega^\bullet (M) \\
& \phi&\mapsto &\phi\circ\omega.
\end{array}
\] (Note the abuse of notation in writing $\omega$ for $\wedge^i\omega$.)
Because $\omega$ satisfies the Maurer-Cartan equation, it is
straightforward to check that $\chi_\omega$ commutes with
differentials, therefore we have a map on the cohomologies:
\[
\chi_\omega: H^\bullet(\frakg)\rightarrow H^\bullet(M).
\]
The cohomology classes in the image of $\chi_\omega$ can be
considered as characteristic classes of the $\frakg$-structure
$\omega$.

In the case of a not necessarily trivialized $\calp$ over $M$,
with the flat connection given by a $G$-equivariant
$\omega\in\Omega^1(\calp)\otimes\frakg$ satisfying the
Maurer-Cartan equation, we analogously obtain a map $\chi_\omega$
\[
\begin{array}{cccc}
\chi_\omega:& \wedge^\bullet \frakg^* &\rightarrow &\Omega^\bullet (\calp) \\
& \phi&\mapsto &\phi\circ\omega.
\end{array}
\]
This again induces a map on cohomologies: $\chi_\omega:
H^\bullet(\frakg)\rightarrow H^\bullet(\calp)$. Should a
trivialization of $\calp$ exist, by choosing one, i.e. picking a
section $s$ of $\calp\rightarrow M$ we obtain a map $s^*\circ
\chi_\omega: H^\bullet(\frakg)\rightarrow H^\bullet(M)$ that
recovers the above case; it depends on $s$ up to homotopy.

In order to deal with the case when a trivialization does not
exist we can also consider the relative version of the above
construction.  Suppose that $K\subset G$ is a subgroup, denote by
$p:\calp\rightarrow\calp/K$ the projection map. The form $\omega$
in $\Omega^1(\calp)\otimes\frakg$ specifying the flat connection
on $\calp$ can be considered as a $G$-equivariant map
$\frakg^*\rightarrow\Omega^1_\calp$.  This restricts to a
$K$-equivariant map $(\frakg/\frakk)^*\rightarrow
p^*\Omega^1_{\calp/K}$.  Thus we have
$(\wedge^\bullet(\frakg/\frakk)^*)^K\rightarrow
(p^*\Omega^\bullet_{\calp/K})^K=\Omega^\bullet_{\calp/K}$. In this
way we obtain the relative version of characteristic classes:
$$\chi_{\omega_K}:H^\bullet(\frakg,K)\rightarrow H^\bullet(\calp/K).$$  As
before, we may pick a section $s$ of $\pi:\calp/K\rightarrow M$
(if it exists) to get the relative characteristic classes map
$s^*\circ \chi_{\omega_K}$ that again depends on the section up to
homotopy. However if $G/K$ is contractible we obtain a canonical
map $(\pi^*)^{-1}\circ \chi_{\omega_K}$.

In order to deal with orbifolds we must adapt the above to the
equivariant setting.  That is we consider a principal $G$-bundle
$\calp$ over $M$ with a $\Gamma$ action.  This means that in
addition to the usual structure we have the action of $\Gamma$ on
$\calp$ and $M$ that makes the projection map
$\Gamma$-equivariant, and the actions of $\Gamma$ and $G$ commute,
i.e. $\gamma(xg)=\gamma(x)g$ for $\gamma\in\Gamma$, $x\in\calp$
and $g\in G$.  If $\calp$ is equipped with a flat connection, we
require the connection to be $\Gamma$-equivariant, i.e. the
connection form $\omega$ is in
$(\Omega^1_\calp\otimes\frakg)^\Gamma$.  We may now proceed as
before with the additional observation that $\chi_{\omega_K}$ is
$\Gamma$-equivariant and thus
$$\chi_{\omega_K}:H^\bullet(\frakg,K)\rightarrow
H^\bullet(\calp/K)^\Gamma$$ and in the case that $G/K$ is
contractible the image of $(\pi^*)^{-1}\circ \chi_{\omega_K}$ is
contained in $H^\bullet(M)^\Gamma$.

Consider a $\gamma\in\Gamma$ and let $M^\gamma=\coprod
M^\gamma_i$, where $M^\gamma_i$ is a $\Gamma^\gamma$-orbit of a
connected component of $M^\gamma$. To any $x\in M^\gamma_i$ we may
associate $\left<\gamma_x\right>$ a conjugacy class in $G$ of
$\gamma_x\in G$ that satisfies $y\gamma_x=\gamma(y)$ for some
choice of $y\in\calp|_x$.  While $\gamma_x$ depends on $y$, its
conjugacy class in $G$ does not. If $\calp$ is equipped with a
$\Gamma$-equivariant flat connection then $\left<\gamma_x\right>$
will be constant on each connected component of $M^\gamma_i$ and
thus on $M^\gamma_i$ itself. Thus for a choice of
$\gamma_i\in\left<\gamma_x\right>$ (notice that $\gamma_i\in G$
unlike $\gamma\in\Gamma$) we may reduce the structure of
$\calp|_{M^\gamma_i}$ to $G^{\gamma_i}$. Explicitly, we consider
the principal $G^{\gamma_i}$-bundle with the fiber over $x\in
M^\gamma_i$ consisting of all $y\in\calp|_x$ such that
$\gamma(y)=y\gamma_i$. Let us denote this principal
$G^{\gamma_i}$-bundle over $M^\gamma_i$ by $\calp^{\gamma_i}$.  It
is evidently $\Gamma^\gamma$-equivariant and it is easy to see
that the $\Gamma$-equivariant connection on $\calp$ will induce a
$\Gamma^\gamma$-equivariant connection on $\calp^{\gamma_i}$.

The characteristic classes formalism above can then be used to
supply the classes associated with $\calp$.  Namely, for every
$\gamma\in\Gamma$ (we need only choose one per conjugacy class),
and a choice of $M^\gamma_i\subset M^\gamma$, choose a
$\gamma_i\in\left<\gamma_x\right>\subset G$ and suppose further
that we have a $K_{\gamma_i}\subset G^{\gamma_i}$ with
$G^{\gamma_i}/K_{\gamma_i}$ contractible. Then we get a
characteristic classes map
$$H^\bullet(\frakg^{\gamma_i},K_{\gamma_i})\rightarrow
H^\bullet(M^\gamma_i)^{\Gamma^\gamma}$$ that essentially does not
depend on the choices of conjugacy class representatives.  We note
that the classes obtained in this way are located inside the
cohomology of the inertia orbifold $\widetilde{X}$.

Now we can apply all this to obtain characteristic classes of
foliations. The key to the application is the consideration of an
appropriate principal $G$-bundle with flat connection.  The
``group" $G$ will be the cross product of the group of formal
coordinates around $0\in\reals^n$ and the group of formal
translations in $\reals^n$. Denote this group by $Diff_n$ and note
that its Lie algebra is $W_n$. Observe that $GL_n(\reals)\subset
Diff_n$ and $Diff_n/GL_n(\reals)$ is contractible.

First, let us recall the non-equivariant situation.  Let $\calf$
be a foliation on a manifold $M$ of codimension $n$. We define the
principal $Diff_n$-bundle $\calp_{Diff}$ over $M$ as follows.  Let
$\hat{\reals}^n_0$ denote the formal neighborhood of $0$ in
$\reals^n$ and $Aut_n$ the group of diffeomorphisms of
$\hat{\reals}^n_0$. With this notation
$Diff_n=Aut_n\ltimes\hat{\reals}^n_0$. Define the $Aut_n$-bundle
$\calp_{Aut}$ by setting the fiber over $x\in M$ to be the set of
diffeomorphisms from $\hat{\reals}^n_0$ to the formal neighborhood
of $x$ in the local leaf space around $x$, i.e. the quotient of
the formal neighborhood of $x\in M$ by the connected components of
the leaves.  Let $\calp_{Diff}=\calp_{Aut}\times_{Aut_n} Diff_n$.
The resulting $Diff_n$-bundle has a canonical flat connection. By
the general theory we get a map:
$$H^\bullet(W_n,GL_n(\reals))\rightarrow H^\bullet(M)$$ that encodes
information about the foliation $\calf$.  More precisely, these
are the real Pontryagin classes of the normal bundle $N\calf$ to
$\calf$ (see Remark \ref{rmk:relgl}).  As pointed out in
\cite{fuchs:book} this gives the well known vanishing theorem of
Bott: monomials of the real Pontryagin classes of $N\calf$ that
have degree $>2n$ are equal to $0$. However one does obtain extra
information using this method, namely since
$GL_n(\reals)/O_n(\reals)$ is contractible we have:
$$H^\bullet(W_n,O_n(\reals))\rightarrow H^\bullet(M)$$ that
extends the map above and thus defines extra classes called
secondary characteristic classes for a general foliation.

We can refine this type of analysis in the special case of a
framed foliation.  Suppose that $N\calf$ is trivialized.  Choosing
a metric yields an exponential map that gives a reduction of
structure of $\calp_{Diff}$ to $F(N\calf)$, the frame bundle of
$N\calf$.  At the same time the trivialization provides a smooth
section of $F(N\calf)$ and since any two metrics are homotopic, we
have a canonical map:
$$H^\bullet(W_n)\rightarrow H^\bullet(M)$$ encoding the data
of the framing as well.

Now we want to use the above idea to define characteristic classes
for a foliation on an orbifold. For simplicity of exposition let
us assume that we are in the setting of a manifold $M$ with a
finite group $\Gamma$ action that preserves a foliation  $\calf$
(of codimension $n$) on $M$.  Then the principal $Diff_n$-bundle
$\calp_{Diff}$ as defined above is automatically
$\Gamma$-equivariant with an equivariant flat connection.  One can
choose a $\Gamma$-equivariant metric on $N\calf$ (any two such
metrics are homotopic), then the exponential map gives a
$\Gamma$-equivariant reduction of structure of $\calp_{Diff}$ to
$F(N\calf)$. This reduction allows us to choose $\gamma_i\in
GL_n(\reals)$. Since $Diff^{\gamma_i}_n/GL^{\gamma_i}_n(\reals)$
is contractible, we have the characteristic classes map:
$$H^\bullet(W_n^{\gamma_i},GL_n^{\gamma_i}(\reals))\rightarrow
H^\bullet(M_0^\gamma)^{\Gamma^\gamma}$$ that can, similarly to the
non-equivariant discussion, be interpreted as a mixture of
Pontryagin and Chern classes for the direct sum decomposition of
the restriction of $N\calf$ to $M^\gamma_i$ induced by the action
of $\gamma$. \emph{We also get an equivariant version of the
vanishing theorem that is considerably stronger.}

To obtain extra information beyond the standard characteristic
classes, observe that
$GL_n^{\gamma_i}(\reals)/O_n^{\gamma_i}(\reals)$ is contractible,
thus we have:
$$\chi_{\gamma,i}: H^\bullet(W_n^{\gamma_i},O_n^{\gamma_i}(\reals))\rightarrow
H^\bullet(M^\gamma_i)^{\Gamma^\gamma}.$$  Explicitly, using the
conventions of Section \ref{sec:realcase}, the action of
$\gamma_i$ on $\reals^n$ gives a decomposition $\reals^n=V_0\oplus
m_{-1}W_{-1}\oplus\bigoplus_{\alpha=1}^k m_\alpha W_\alpha$.  Then
$GL_n^{\gamma_i}(\reals)=GL(V_0)\times
GL_{m_{-1}}(\reals)\times\prod GL_{m_\alpha}(\C)$ and
$O_n^{\gamma_i}(\reals)=O(V_0)\times O_{m_{-1}}(\reals)\times\prod
U_{m_\alpha}$.

\begin{dfn}
For a foliation $\calf$ on an orbifold $X$, we call elements in
the image of the map
$$\chi=\oplus
\chi_{\gamma, i}: \bigoplus_{\left<\gamma\right>\subset\Gamma,\,
i}H^\bullet(W_n^{\gamma_i},O_n^{\gamma_i}(\reals))\rightarrow
H^\bullet(\tilde{X})$$ characteristic classes of $\calf$.
\end{dfn}

The cohomology of
$H^\bullet(W_n^{\gamma_i},O_n^{\gamma_i}(\reals))$  is computed in
Theorem \ref{thm:top-real}. It is isomorphic to the cohomology of
the topological space $X_\rho^{\gamma_i}/O_n^{\gamma_i}$, which is
fibered over $B_\rho^{\gamma_i}$. The cohomology ring of
$B_\rho^{\gamma_i}$ is the polynomial ring with
$dim(V_0)+m_{-1}+2\sum m_\alpha$ generators truncated at degree
$>2dim(V_0)$. The  images under $\chi$ of these generators
correspond to the Pontryagin and Chern classes of the
decomposition of the restriction of the normal bundle $N\calf$.
Furthermore, we notice that
$X_\rho^{\gamma_i}/O_n^{\gamma_i}\rightarrow B_\rho^{\gamma_i}$
has a cohomologically nontrivial fiber $U_{dim(V_0)}\times
U_{m_{-1}}\times \prod(U_{m_\alpha}\times U_{m_\alpha})/
O(V_0)\times O_{m_{-1}}(\reals)\times\prod U_{m_\alpha}$. Thus the
cohomology group of $X_\rho^{\gamma_i}/O_n^{\gamma_i}$ is larger
than the one of $B_\rho^{\gamma_i}$. The extra cohomology classes
give rise to the secondary characteristic classes of $\calf$. For
example, looking at the $E_2$ term of the spectral sequence
associated to the fibration, we see that any class in the upper
right hand corner will survive as a class in the cohomology of
$X_\rho^{\gamma_i}/O_n^{\gamma_i}$.  These can be thought of, in a
sense, as generalizations of the Godbillon-Vey class.

\begin{rmk}
Different approaches such as \cite{kt:foliation} can also be used
to define characteristic classes for foliations on orbifolds.
Here, the Gelfand-Fuchs cohomology approach is used as an
application of the Lie algebra cohomology computation.
\end{rmk}

\vskip 1cm \noindent Department of Mathematics, University of
California, Davis, CA, USA \newline \emph{E-mail address}:
\textbf{ishapiro@math.ucdavis.edu}
\newline \emph{E-mail address}:
\textbf{xtang@math.ucdavis.edu}

\end{document}